\documentclass[letterpaper, 10 pt, conference]{ieeeconf}  

\IEEEoverridecommandlockouts                              

\usepackage{graphics} 
\usepackage{epsfig} 
\usepackage{graphicx} 
\usepackage{times} 
\usepackage{amsmath} 
\usepackage{amssymb}  
\usepackage{algorithm,algorithmic}
\usepackage{subfigure}

\title{\LARGE \bf
Initial Error Tolerant Distributed Mean Field Control under Partial and Discrete Information 
}

\author{Yuxin Jin, Haotian Wang, Wang Yao and Xiao Zhang 
\thanks{* This work has been submitted to the IEEE for possible publication. Copyright may be transferred without notice, after which this version may no longer be accessible.}
\thanks{* This work was supported by the National Science and Technology Major Project (Grant No. 2022ZD0116401) and the Research Funding of Hangzhou International Innovation Institute of Beihang University (Grant No. 2024KQ161). (Corresponding author: Wang Yao and Xiao Zhang).}
\thanks{Y. Jin is with ShenYuan Honors College and School of Mathematical Sciences, Beihang University, Beijing 100191, China; Key Laboratory of Mathematics, Informatics and Behavioral Semantics, Ministry of Education, Beihang University, Beijing 100191, China (e-mail: yxjin@buaa.edu.cn).}
\thanks{H. Wang is with School of Mathematical Sciences, Beihang University, Beijing 100191, China (e-mail: williamelwht@buaa.edu.cn).}%
\thanks{ W. Yao is with School of Artificial Intelligence and LMIB, Beihang University, Beijing 100191, China; Hangzhou International Innovation Institute of Beihang University, Hangzhou 311115, China (e-mail: yaowang@buaa.edu.cn).}%
\thanks{  X. Zhang is with the School of Mathematical Sciences, Beihang University, Beijing 100191, China; Key Laboratory of Mathematics, Informatics and Behavioral Semantics, Ministry of Education, Beihang University, Beijing 100191, China (e-mail: xiao.zh@buaa.edu.cn).}%
}

\begin{document}

\maketitle
\thispagestyle{empty}
\pagestyle{empty}

\begin{abstract}
In this paper, an initial error tolerant distributed mean field control method under partial and discrete information is introduced, where each agent only has discrete observations on its own state. First, we study agents' behavior in linear quadratic mean field games (LQMFGs) under heterogeneous erroneous information of the initial mean field state (MF-S), and formulate the relationships between initial errors and systemic deviations. Next, by capturing the initial error affection on the private trajectory of an agent, we give a distributed error estimation method based on maximum likelihood estimation (MLE), where each agent estimates information errors only based on discrete observations on its private trajectory. Furthermore, we establish an error-based segmented state estimation method, design the initial error tolerant distributed mean field control method (IET-DMFC), and demonstrate the consistent property of state estimation as observation frequency increases. Finally, simulations are performed to verify the efficiency of the algorithm and the consistent properties. 
\end{abstract}

\section{Introduction}

Mean field games (MFGs), independently introduced by Huang-Caines-Malhamé \cite{1}\cite{2} and Lasry-Lions \cite{3}-\cite{5}, is a pivotal framework for modeling large-scale systems, where interactions among agents are simplified through statistical aggregation into mean field (MF) terms. In an infinite population, MF control laws admit each agent only using its private states and parameters of the system to predict MF and give its optimal control, enabling decentralized strategies which yield Nash equilibrium. This paradigm has been successfully applied across diverse domains, such as smart grids, autonomous systems and crowd motion \cite{15}-\cite{18}.  

MF control laws require agents to possess accurate initial MF information to derive optimal strategies. However, this requirement often fails in practice. Observation errors, environmental noise, and sensor interference \cite{19}\cite{20} can lead to biased initial estimates, causing non-optimality. Furthermore, communication constraints are commonly applied in multi-agent systems \cite{21}\cite{22}, and large population scales causes difficulties on observing MF, leading to partial information, making it difficult for agents to correct information errors through communications and process observations on MF. 

Several works have studied MFGs and MF control laws under partial or erroneous information. In \cite{6}-\cite{11}, partial observation situation is discussed, where \cite{6} and \cite{11} analyze agents' partial observation on their own states, \cite{12} focuses on discrete observations on MF, and \cite{7}-\cite{10} consider a partially observed major player in major-minor MFGs. For erroneous information, paper \cite{13} investigates erroneous initial information affection and error correction method in multi-population MFGs, and \cite{14} proposes MF stochastic control, in which each agent observes a random subset of agents, estimates own and population dynamical parameters.  

\begin{figure}[h!]
\centering
\includegraphics[width=3in]{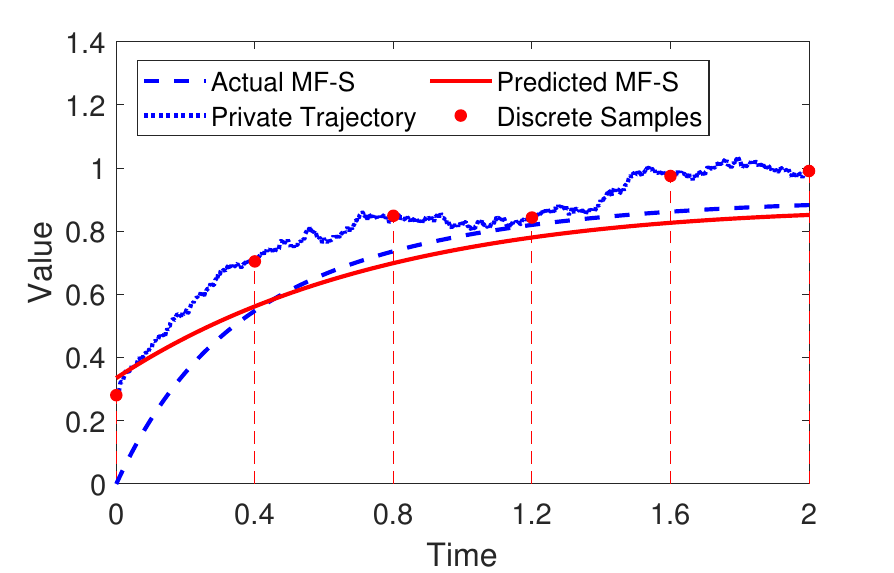}
\caption{Schematic of an one-dimensional scene under partial and discrete information. An selected agent has access to discrete samples of its own private trajectory and predicted MF-S that are marked by red. Actual MF-S and continuous private trajectory marked by blue are not accessible. }
\label{fig_0}
\end{figure}

Despite these advancements, approaches in MF control laws often rely on process observations on other agents, or lack error-tolerant mechanisms under stochastic dynamics. To address these challenges while adhering to practical constraints of discrete sampling,  we propose an initial error tolerant distributed mean field control (IET-DMFC) framework. As depicted in Fig.\ref{fig_0}, each agent only has access to strictly local information: discrete observations of its private states and potentially erroneous initial information. At the initial moment, agents have heterogeneous erroneous observations on initial MF and compute their feedback control laws. The feedback control laws are modified at predefined time series. At each time point, each agent utilizes its time-equidistant recorded private trajectory observations to estimate the current MF-S and refine its control strategies accordingly.

The main contributions of this paper can be listed as follows:
\begin{itemize}
  \item Based on the analysis on heterogeneous initial error affections on agents' private trajectory, we design a distributed error estimation method based on MLE where agents infer errors solely from discrete observations of their private trajectories. 
  \item Applying the relationships between estimation errors and population evolution, we design an error-based segmented state estimation method based on the distributed error estimation method.
  \item We establish the IET-DMFC framework. The consistent property of state estimations related to observation frequency is demonstrated.  
\end{itemize}

The rest of paper is organized as follows: in Section 2, we discuss agents' behavior and error propagation in LQMFGs under heterogeneous erroneous information. In Section 3, we propose a distributed error estimation method to correct information errors. In Section 4, we propose IET-DMFC method, and discuss its consistent property.  In Section 5, we conduct simulations and verify our conclusions. 

\section{Mean Field Games under Heterogeneous Erroneous Information}
In this section, we introduce the LQMFG model where agents get heterogeneous observations on initial MF-S $z_0$. Agents give their feedback control law at the initial moment. We consider a game with $N$ agents, $\mathcal{A}_i,1\leq i\leq N$, $N\rightarrow \infty$, the parameter set of the system is $\Theta$. Select $\mathcal{A}_i$ as a generic agent.
\subsection{Assumptions and Scenario Description}
\emph{A1}: 
Initial private state $x_i(0)$, system parameter set $\Theta$ are accessible for $\mathcal{A}_{i}$.  

\emph{A2}: 
Initial observation $z_i(0)=z^0+E_i$ is the $z_0$ as observed by $\mathcal{A}_{i}$ with private error $E_i \in \mathbb{R}^n$. Define mean error as $\bar{E}:=\frac{1}{N}\Sigma_{i=1}^NE_i$, $(\bar{E},E_i)\in \Lambda$, the error set $\Lambda\subset\mathbb{R}^{2n}$ is an open and compact set, $0\in \Lambda$.

\emph{A3}: 
$\mathcal{A}_i$ takes $z_i(0)$ as the correct initial MF-S, and assumes all agents have correct observations on initial MF-S to give its strategy.

\emph{A4}: 
At $t=0$, $\mathcal{A}_{i}$ gives its feedback control law $u_{i}(t)=\phi_i(x_i(t),t)$ and evolves according to this strategy.

\subsection{Dynamics and Cost Functions}
Let $(\Omega, \mathcal{F},\mathbb{P})$ be a complete probability space and $T>0$. Suppose that $W_i,1\leq i \leq N$ are independent $n$-dimensional standard Wiener processes defined on $(\Omega, \mathcal{F},\mathbb{P})$, and $x_0^i, 1\leq i \leq N$ are $n$-dimensional vectors. Let $(x_i(t))_{0\leq t\leq T}$ be the state of $\mathcal{A}_i$, and $x_0^i$ represents the initial state of $\mathcal{A}_i$. 

The dynamics of $\mathcal{A}_i$ are given by
\begin{equation}
\label{deqn_ex1}
\begin{split}
dx_i(t)=&[Ax_i(t)+Bu_i(t)+Cz(t)+F\bar{u}(t)]dt+DdW_i(t),\\
x_i(0)=&x_0^i,
\end{split}
\end{equation}
where $A, B, C, F, D$ are matrices of suitable sizes, $z(t):=1/N\Sigma_{j=1}^N x_j(t)$ and $\bar{u}(t):=1/N\Sigma_{j=1}^N u_j(t)$ are the MF-S and mean field control (MF-C). Control $u_i(t)$ is in $L^2_{\mathcal{F}_i}(0,T;\mathbb{R}^m)$, $\mathcal{F}^i_t:=\sigma(x_i(0),W_i(s), s\leq t)$.

The cost functional of $\mathcal{A}_i$ is 
\begin{equation}
\label{deqn_ex2}
\begin{split}
J(u_i)=&\frac{1}{2}\mathbb{E}[\int_{0}^{T}[\|x_i(t)-s\|_{Q_{I}}^2+\|u_i(t)\|_R^2+\|x_i(t)-(\Gamma\\
&z(t)+\eta)\|_Q^2]dt+\|x_i(T)-\bar{s}\|_{\bar{Q}_{I}}^2+\|x_{i}(T)-(\bar{\Gamma}\\
&z(T)+\bar{\eta})\|_{\bar{Q}}^2],
\end{split}
\end{equation}
where we define $\|X\|_Q^2=X^TQX$. $Q_{I}, Q, R, \bar{Q}_{I}, \bar{Q}$ are square matrices, and $R$ is positive definite. Let $\Theta$ be the set of parameters of the above systems.
\subsection{Optimal Control}
As $N\rightarrow \infty$, for given continuous deterministic process $(z_i(t))_{0\leq t\leq T}$ and $(\bar{u}_i(t))_{0\leq t\leq T}$, according to the stochastic maximum principle \cite{23}\cite{25}, $\mathcal{A}_i$ can give its optimal control as follows.

\textbf{Lemma 2.1} $\mathcal{A}_i$'s optimal control $u_i(t)=-R^{-1}B^Tp_i(t)$, where \\
\begin{equation}
\label{deqn_ex3}
\begin{split}
dx_i(t)=&(Ax_i(t)-BR^{-1}B^Tp_i(t)+Cz_i(t)+F\bar{u}_i(t))dt\\
&+DdW_i(t)\\
x_i(0)=&x_0^i\\
dy_i(t)=&(-A^Ty_i(t)+Q\Gamma z_i(t)-(Q_I+Q)x_i(t)+Q_Is\\
&+Q\eta)dt\\
y_i(T)=&(\bar{Q}_I+\bar{Q})x_i(T)-\bar{Q}_I\bar{s}-\bar{Q}(\bar{\Gamma}z_i(T)+\bar{\eta}).
\end{split}
\end{equation}
such that $p_i(t)=\mathbb{E}[y_i(t)|\mathcal{F}_t^i]$.

\textbf{Remark 2.1} The optimal control $u_i$ has a feedback representation $-R^{-1}B^T(P_1(t)x_i(t)+g_i(t))$, where $P_1(t)$ satisfies a non-symmetric riccati equation
\begin{equation}
\label{deqn_ex4}
\begin{split}
-dP_1(t)=&[P_1A+A^TP_1+(Q_I+Q)-P_1BR^{-1}B^TP_1]dt,\\
P_1(T)=&\bar{Q}_I+\bar{Q},
\end{split}
\end{equation}
and $g_i(t)$ satisfies 
\begin{equation}
\label{deqn_ex5}
\begin{split}
dg_i(t)=&-[(A^T-P_1BR^{-1}B^T)g_i(t)+(P_1C-Q\Gamma)z_i(t)+\\
&P_1F\bar{u}_i(t)-Q_Is-Q\eta]dt,\\
g_i(T)=&-\bar{Q}_I\bar{s}-\bar{Q}(\bar{\Gamma}z_i(T)+\bar{\eta}).
\end{split}
\end{equation}

\subsection{Agents' Behavior under Complete Information}
We first consider the complete information situation, where all agents give the consistent correct predictions $(z(t))_{0\leq t\leq T}$ and $(\bar{u}(t))_{0\leq t\leq T}$ for MF. 

On the one hand, a Nash equilibrium is reached if and only if each agent's control is the optimal response to the current mean field term. On the other hand, $z(t)$ and $\bar{u}(t)$ are generated by agents' behavior. According to the last subsection, we have the following theorem:

{\bf{Theorem 2.2}} MF-S $z(t)$ and MF-C $\bar{u}(t)=-R^{-1}B^Tp(t)$ under complete information satisfy the following equations
\begin{equation}
\label{deqn_ex6}
\begin{split}
&d\begin{pmatrix} 
z\\
p\\
\end{pmatrix}=\left\{\begin{pmatrix} 
A+C & -(B+F)R^{-1}B^T \\
\mathcal{Q} & -A^T
\end{pmatrix}
\begin{pmatrix} 
z\\
p\\
\end{pmatrix}-
\begin{pmatrix} 
0\\
\nu\\
\end{pmatrix}\right\}dt,\\
&z(0)=z^0,\\
&p(T)=(\bar{Q}_I+\bar{Q}-\bar{Q}\bar{\Gamma})z(T)-\bar{Q_I}\bar{s}-\bar{Q}\bar{\eta}.\\
\end{split}
\end{equation} 
where $\nu=Q_Is+Q\eta,\mathcal{Q}=Q\Gamma-Q_I-Q$.

{\bf{Remark 2.2}} We notice that $p(t)=P_0(t)z(t)+\mathcal{G}(t)$, where $P_0(t)$ satisfies a non-symmetric riccati equation
\begin{equation}
\label{deqn_ex7}
\begin{split}
-dP_0(t)=&\{P_0(t)(A+C)+A^TP_0(t)+(Q_I+Q-Q\Gamma)\\
&-P_0(t)(B+F)R^{-1}B^TP_0(t)\}dt,\\
P_0(T)=&\bar{Q}_I+\bar{Q}-\bar{Q}\bar{\Gamma}.\\
\end{split}
\end{equation}
and $\mathcal{G}(t)$ satisfies the backward ordinary differential equation
\begin{equation}
\label{deqn_ex8}
\begin{split}
d\mathcal{G}(t)=&\{-(A^T-P_0(t)(B+F)R^{-1}B^T)\mathcal{G}+Q_Is+Q\eta\}dt,\\
\mathcal{G}(T)=&-\bar{Q_I}\bar{s}-\bar{Q}\bar{\eta}.\\
\end{split}
\end{equation}
According to Remark 2.1, $p(t)=P_1(t)z(t)+g(t)$.

When (\ref{deqn_ex7}) and (\ref{deqn_ex4}) have unique solutions, $\mathcal{A}_i$ can get its feedback control law. At the initial time, $P_0$, $P_1$ and $\mathcal{G}$ can be computed by $\mathcal{A}_i$, $\forall i$. For given $z^0$, agents can compute the mean field state(MF-S) $z(t)$ and mean field control(MF-C) $\bar{u}(t)$ for $0\leq t\leq T$. 

Substitute $P_0(t)$ and $\mathcal{G}(t)$ into (\ref{deqn_ex6}), then 
\begin{equation}
\label{deqn_ex9}
\begin{split}
dz(t)=&[(A+C-(B+F)R^{-1}B^TP_0(t))z(t)-(B+F)\\
&R^{-1}B^T\mathcal{G}(t)]dt,\\
z(0)=&z^0,\\
\bar{u}(t)=&-R^{-1}B^T(P_0(t)z(t)+\mathcal{G}(t))
\end{split}
\end{equation} 

For $(z(t))_{0\leq t\leq T}$, $\mathcal{A}_i$ can solve (\ref{deqn_ex5}) for $(g(t))_{0\leq t\leq T}$. According to Remark 2.1, $\mathcal{A}_i$'s feedback optimal control is
\begin{equation}
\label{deqn_ex10}
\begin{split}
&u_i(t)=\phi(x_i(t),t), 0\leq t\leq T\\
&\phi(x_i(t),t)=-R^{-1}B^T(P_1(t)x_i(t)+g(t)).\\
\end{split}
\end{equation} 

Then $\mathcal{A}_i$'s behavior can be described by Algorithm 1.
\begin{algorithm}[!ht]
	\caption{$\mathcal{A}_i$'s Behavior under Complete Information}
    \label{power}
    \begin{algorithmic}[1] 
        \STATE Solve (\ref{deqn_ex7}), (\ref{deqn_ex8});
        \STATE Predict $(z(t))_{0\leq t\leq T}$ by (\ref{deqn_ex9});
        \STATE Get $(g(t))_{0\leq t\leq T}$ by (\ref{deqn_ex5}) ;
        \STATE $u_i(t)\leftarrow\phi(x_i(t),t), 0\leq t\leq T$;
        \STATE Evolve according to (\ref{deqn_ex1});
    \end{algorithmic}
\end{algorithm}

\subsection{Agents' Behavior under Erroneous Information}
Consider the situation where agents have heterogeneous erroneous information of initial MF-S.
According to \emph{A3}, MF-S predicted by $\mathcal{A}_i$ is the mean field equilibrium under complete information where $z(0)=z_i(0)$, which satisfies
\begin{equation}
\label{deqn_ex11}
\begin{split}
dz_i(t)=&[(A+C-(B+F)R^{-1}B^TP_0(t))z_i(t)-(B+F)\\
&R^{-1}B^T\mathcal{G}(t)]dt,\\
z_i(0)=&z^0+E_i,\\
\bar{u}_i(t)=&-R^{-1}B^T(P_0(t)z_i(t)+\mathcal{G}(t)).
\end{split}
\end{equation} 

The initial private error $E_i$ affects predicted MF-S $(z_i(t))_{0\leq t\leq T}$, and then affects $(g_i(t))_{0\leq t\leq T}$ through (\ref{deqn_ex5}). Then $\mathcal{A}_i$'s behavior can be described by Algorithm 2.
\begin{algorithm}[!ht]
	\caption{$\mathcal{A}_i$'s Behavior under Erroneous Information}
    \label{power}
    \begin{algorithmic}[1] 
        \STATE Solve (\ref{deqn_ex7}), (\ref{deqn_ex8});
        \STATE Predict $(z_i(t))_{0\leq t\leq T}$ by (\ref{deqn_ex11});
        \STATE Get $(g_i(t))_{0\leq t\leq T}$ by (\ref{deqn_ex5}) ;
        \STATE $u_i(t)\leftarrow\phi_i(x_i(t),t), 0\leq t\leq T$;
        \STATE Evolve according to (\ref{deqn_ex1});
    \end{algorithmic}
\end{algorithm}
\subsection{Actual MF}
Given agents' behavior under erroneous information, the actual evolution of the game can be described. Define $(z_A(t))_{0\leq t\leq T}$ as the actual mean field state under erroneous information, according to agents' dynamics (\ref{deqn_ex1}), when $N\rightarrow\infty$, we have
\begin{equation}
\label{deqn_ex12}
\begin{split}
dz_A(t)=&[(A+C-(B+F)R^{-1}B^TP_1(t))z_A(t)-(B+\\
&F)R^{-1}B^T\bar{g}(t)]dt,\\
d\bar{z}(t)=&[(A+C-(B+F)R^{-1}B^TP_0(t))\bar{z}(t)-(B+F)\\
&R^{-1}B^T\mathcal{G}(t)]dt,\\
d\bar{g}(t)=&-[(A^T-P_1(t)BR^{-1}B^T)\bar{g}(t)+(P_1C-Q\Gamma)\\
&\bar{z}(t)-P_1FR^{-1}B^T(P_1(t)\bar{z}(t)+\bar{g}(t))-Q_Is\\
&-Q\eta]dt,\\
z_A(0)=&z^0,\\
\bar{z}(0)=&z^0+\bar{E},\\
\bar{g}(T)=&-\bar{Q}_I\bar{s}-\bar{Q}(\bar{\Gamma}\bar{z}(T)+\bar{\eta}),\\
\end{split}
\end{equation} 
where $\bar{g}(t):=\frac{1}{N}\sum_{i=1}^Ng_i(t)$, and $\bar{z}(t):=\frac{1}{N}\sum_{i=1}^Nz_i(t)$.

\section{Distributed Error Estimation with Partial and Discrete Information}
In this section, for the scene in the last section, we give a distributed error identification method where $\mathcal{A}_i$ only has access to discrete observations on its private trajectory $x_i$. We set observation time series as $\mathbb{T}:={t_0,t_1,t_2,...,t_k}$, where $t_0=0$. $\mathcal{A}_i$'s information set is $\{x_i(t_j),j=0,1,2,...,k\}\bigcup{z_i(0)}\bigcup\Theta$.

First, we discuss the information error affection on agents' private trajectories. Then, applying the relationships between information error and the deviations of observations from predictions, we turn the error identification problem into a parameter estimation problem. Finally, for this parameter estimation problem, we give the distributed error estimation method based on MLE.     
\subsection{Information Error Affections}
We first consider the information error affections on predicted MF and actual MF. 

Remind $(z(t))_{0\leq t\leq T}$ and $(\bar{u}(t))_{0\leq t\leq T}$ as MF, and $\phi(\cdot,t)=-R^{-1}B^T(P_1(t)\cdot+g(t))$ as the feedback control law under complete information. 

Define $\Delta z_i(t):=z_i(t)-z(t)$, $\Delta \bar{z}(t):=\bar{z}(t)-z(t)$, then according to (\ref{deqn_ex11}) and (\ref{deqn_ex9}), we have
\begin{equation}
\label{deqn_ex13}
\begin{split}
\Delta z_i(t)=&\Phi_1(t)\Phi_1^{-1}(0)E_i,\\
\Delta \bar{z}(t)=&\Phi_1(t)\Phi_1^{-1}(0)\bar{E},\\
\end{split}
\end{equation} 
where $\Phi_1(t)$ is the solution to the following matrix differential equation
\begin{equation}
\label{deqn_ex14}
\begin{split}
d\Phi_1(t)=&[(A+C-(B+F)R^{-1}B^TP_0(t))\Phi_1(t)]dt,\\
\Phi_1(0)=&I.\\
\end{split}
\end{equation} 

Then for $\Delta \bar{g}(t):=\bar{g}(t)-g(t)$, $\Delta g_i(t):=g_i(t)-g(t)$, according to (\ref{deqn_ex5}), using the method of variation of constant, we have
\begin{equation}
\label{deqn_ex15}
\begin{split}
\Delta g_i(t)=&\mathcal{M}_g(t)E_i,\\
\Delta \bar{g}(t)=&\mathcal{M}_g(t)\bar{E},\\
\end{split}
\end{equation} 
where
\begin{equation}
\label{deqn_ex16}
\begin{split}
\mathcal{M}_g(t)=&-\Phi_g(t)[\Phi_g^{-1}(T)\bar{Q}\bar{\Gamma}\Phi_1(T)\Phi_1^{-1}(0)+\int_T^t\Phi_g^{-1}(s)\\
&(P_1C-P_1FR^{-1}B^TP_1-Q\Gamma)\Phi_1(s)\Phi_1^{-1}(0)ds],\\
\end{split}
\end{equation} 
and $\Phi_g(t)$ is the solution to the following matrix differential equation
\begin{equation}
\label{deqn_ex17}
\begin{split}
d\Phi_g(t)=&-[(A^T-P_1(t)(B+F)R^{-1}B^T)\Phi_g(t)]dt,\\
\Phi_g(T)=&I.\\
\end{split}
\end{equation} 

Finally, for $\Delta z_A(t):=z_A(t)-z(t)$, according to (\ref{deqn_ex5}), using the method of variation of constant, we have
\begin{equation}
\label{deqn_ex18}
\begin{split}
\Delta z_A(t)=&\mathcal{M}_z(t)\bar{E},\\
\mathcal{M}_z(t)=&-\Phi_z(t)\int_0^t\Phi_z^{-1}(s)(B+F)R^{-1}B^T\mathcal{M}_g(s)ds,\\
\end{split}
\end{equation} 
where $\Phi_g(t)$ is the solution to the following matrix differential equation
\begin{equation}
\label{deqn_ex19}
\begin{split}
d\Phi_z(t)=&[(A+C-(B+F)R^{-1}B^TP_1(t))\Phi_z(t)]dt,\\
\Phi_z(0)=&I.\\
\end{split}
\end{equation} 
\subsection{Transition Probability Density of $\mathcal{A}_i$'s Private Trajectory}
To identify the information error, $\mathcal{A}_i$ can compare its predictions and observations. We go back to (\ref{deqn_ex1}), the actual private trajectory of $\mathcal{A}_i$ satisfies
\begin{equation}
\label{deqn_ex20}
\begin{split}
dx_i(t)=&[(A-BR^{-1}B^TP_1(t))x_i(t)-BR^{-1}B^Tg_i(t)+\\
&(C-FR^{-1}B^TP_1(t))z_A(t)-FR^{-1}B^T\bar{g}(t)]dt\\
&+DdW_i(t),\\
x_i(0)=&x_0^i,
\end{split}
\end{equation}
where $z_A(t)$ and $\bar{g}(t)$ are unknown. But according to last subsection, we can replace them by information error and computable variables 
\begin{equation}
\label{deqn_ex21}
\begin{split}
z_A(t)=&z_i(t)+\mathcal{M}_z(t)\bar{E}-\Phi_1(t)\Phi_1(0)^{-1}E_i,\\
\bar{g}(t)=&g_i(t)+\mathcal{M}_g(t)(\bar{E}-E_i).\\
\end{split}
\end{equation}

Then $\mathcal{A}_i$'s dynamics can be rewritten as
\begin{equation}
\label{deqn_ex22}
\begin{split}
dx_i(t)=&[(A-BR^{-1}B^TP_1(t))x_i(t)-(B+F)R^{-1}B^T\\
&g_i(t)+(C-FR^{-1}B^TP_1(t))z_i(t)]dt\\
&+[\bar{\mathcal{K}}(t)\bar{E}+\mathcal{K}_i(t)E_i]dt+DdW_i(t),\\
x_i(0)=&x_0^i.
\end{split}
\end{equation}
where 
\begin{equation}
\label{deqn_ex23}
\begin{split}
\bar{\mathcal{K}}(t)=&(C-FR^{-1}B^TP_1(t))\mathcal{M}_z(t)-FR^{-1}B^T\mathcal{M}_g(t),\\
\mathcal{K}_i(t)=&FR^{-1}B^T\mathcal{M}_g(t)-(C-FR^{-1}B^TP_1)\Phi_1(t)\Phi_1^{-1}(0)\\
\end{split}
\end{equation}

Then we have the following theorem which can help us establish our error estimation method. Define $E=(\bar{E}^T,E_i^T)^T$, and $\mathcal{K}(t)=(\bar{\mathcal{K}}(t),\mathcal{K}_i(t))$.

{\bf{Theorem 3.1}} For given $\bar{E},E_i$, the transition probability density of $\mathcal{A}_i$'s private trajectory $x_i(t)$ from $s$ to $t$, $t>s>0$ is
\begin{equation}
\label{deqn_ex24}
\begin{split}
\mathcal{T}(t,y|s,x;E)=\frac{1}{\sqrt{(2\pi)^d|\Sigma|}}\exp^{-\frac{1}{2}(y-\mu)^T\Sigma^{-1}(y-\mu)}
\end{split}
\end{equation}
which is a multivariate normal distribution with mean $\mu(t|s,x;E)$ and covariance $\Sigma(t|s)$:
\begin{equation}
\label{deqn_ex25}
\begin{split}
&\mu(t|s,x;E)=F(t|s,x)+\Phi_x(t)\int_s^t\Phi_x^{-1}(\omega)\mathcal{K}(\omega)d\omega E,\\
&\Sigma(t|s)=\Phi_x(t)\int_s^t\Phi_x^{-1}(\omega)DD^T\Phi_x^{-T}(\omega)d\omega \Phi_x^T(t),
\end{split}
\end{equation}
where
\begin{equation}
\label{deqn_ex26}
\begin{split}
&F(t|s,x)=\Phi_x(t)\Phi_x^{-1}(s)x+\Phi_x(t)\int_s^t\Phi_x^{-1}(\omega)\mathcal{B}_i(\omega)d\omega,\\
&\mathcal{B}_i(t)=(C-FR^{-1}B^TP_1)z_i(t)-(B+F)R^{-1}B^Tg_i(t),
\end{split}
\end{equation}
and $\Phi_x(t)$ is the solution to the following matrix differential equation
\begin{equation}
\label{deqn_ex27}
\begin{split}
d\Phi_x(t)=&-[(A-BR^{-1}B^TP_1(t))\Phi_x(t)]dt,\\
\Phi_x(0)=&I.\\
\end{split}
\end{equation}

Proof: 

Applying definitions to (\ref{deqn_ex22}), using the method of variation of constant, we have
\begin{equation*}
\begin{split}
y=&F(t|s,x)+\Phi_x(t)\int_s^t\Phi_x^{-1}(\omega)\mathcal{K}(\omega)d\omega E\\
&+\Phi_x(t)\int_s^t\Phi_x^{-1}(\omega)DdW_i(\omega),
\end{split}
\end{equation*}
where 
\begin{equation*}
\begin{split}
&\Phi_x(t)\int_s^t\Phi_x^{-1}(\omega)DdW_i(\omega)\sim\mathcal{N}_n(0,\Sigma(t|s)),\\
&\Sigma(t|s)=\Phi_x(t)\int_s^t\Phi_x^{-1}(\omega)DD^T\Phi_x^{-T}(\omega)d\omega \Phi_x^T(t).\\
\end{split}
\end{equation*}
Since $F(t|s,x), \Phi_x(t),\mathcal{K}(t)$ are all deterministic, we have
\begin{equation*}
\begin{split}
y\sim\mathcal{N}_n(\mathbb{E}[t,y|s,x;E],\Sigma(t|s)),
\end{split}
\end{equation*}
where
\begin{equation*}
\begin{split}
\mathbb{E}[t,y|s,x;E]=F(t|s,x)+\Phi_x(t)\int_s^t\Phi_x(\omega)\mathcal{K}(\omega)d\omega E.
\end{split}
\end{equation*}
$\Box$

\subsection{Likelihood Function}
Above discussion can turn the initial error estimation problem into a parameter estimation problem of a stochastic differential equation.

{\bf{Problem 3.1}} For known observation value $x_i(t_j),j=0,1,2,...,k$, $\Theta, z_i(0)$, estimate $E$. where $x_i$ satisfies
\begin{equation*}
\begin{split}
dx_i(t)=&[(A-BR^{-1}B^TP_1(t))x_i(t)+\mathcal{B}_i(t)]dt\\
&+\mathcal{K}(t)Edt+DdW_i(t).\\
\end{split}
\end{equation*}

According to theorem 3.1, we can give the likelihood function
\begin{equation}
\label{deqn_ex28}
\begin{split}
\mathcal{L}(E|x_i)=\prod_{j=1}^k\mathcal{T}(t_j,x_i(t_j)|t_{j-1},x_i(t_{j-1});E),
\end{split}
\end{equation}
and the log-likelihood function
\begin{equation}
\label{deqn_ex29}
\begin{split}
l(E|x_i)=\sum_{j=1}^klog(\mathcal{T}(t_j,x_i(t_j)|t_{j-1},x_i(t_{j-1});E)).
\end{split}
\end{equation}

\subsection{Distributed Error Estimation based on MLE}
Define the maximum likelihood estimation for $E$ as $\hat{E}$, then  
\begin{equation}
\label{deqn_ex30}
\begin{split}
\hat{E}=\arg\max_{E\in \Lambda}l(E|x_i).
\end{split}
\end{equation}

Since the form of the likelihood function $\mathcal{L}(E|x_i)$ can be computed only according to $\Theta$ and $z_i(0)$, $\mathcal{A}_i$ can compute this likelihood function and estimate $E$ by solving the above equation for $\hat{E}$. 
The following algorithm shows the error estimation behavior of $\mathcal{A}_i$.

\begin{algorithm}[!ht]
    \renewcommand{\algorithmicrequire}{\textbf{Input: }}
	\renewcommand{\algorithmicensure}{\textbf{Output:}}
	\caption{Error Estimation Method for $\mathcal{A}_i$}
    \label{power}
    \begin{algorithmic}[1] 
        \REQUIRE  $z_i(0),\Theta,x_i(t_j),j=0,...,k$; 
	    \ENSURE $\hat{E}$; 
        \STATE Solve (\ref{deqn_ex7}), (\ref{deqn_ex8}), (\ref{deqn_ex11}) for $z_i,P_0,\mathcal{G}$;
        \STATE Solve (\ref{deqn_ex4}), (\ref{deqn_ex5}) for $P_1,g_i$;
        \STATE Solve (\ref{deqn_ex14}), (\ref{deqn_ex17}), (\ref{deqn_ex19}), (\ref{deqn_ex27}) for $\Phi_1,\Phi_z,\Phi_g,\Phi_x$;
        \STATE Solve (\ref{deqn_ex31}) for $\hat{E}$.
        \STATE \textbf{return} $\hat{E}$.
    \end{algorithmic}
\end{algorithm}

According to the conclusions on consistent properties of MLE discussed in \cite{24} and \cite{14}, we have the following lemma:

{\bf{Lemma 3.1}} When $A1-A4$ holds, $t_j=jT/k$, we have
\begin{equation}
\label{deqn_ex31}
\begin{split}
\hat{E}\stackrel{a.s.}{\longrightarrow} E \ as\ k\rightarrow \infty
\end{split}
\end{equation}

Which is, when observations are time-equidistant, $\hat{E}$ converges to real error value $E$ with probability $1$ as observation frequency tends to infinity.

\section{Initial Error Tolerant Distributed Mean Field Control under Partial and Discrete Information}
According to last section, $\mathcal{A}_i$ can estimate its initial private error $E_i$ and mean error $\bar{E}$ by solving a parameter estimation problem, only using discrete observations of its private trajectory $x_i(t)$. 

However, not only error correction, but also strategy modification is needed to improve the system performance. Therefore we consider a segmented strategy modification scene, where agents modify their strategies at a series of time points, have no process observations on MF, and only record discrete samples of their own private trajectories.

In this section, we establish an error-based state estimation method, and give our initial error tolerant distributed mean field control algorithm, in which each agent estimates the actual MF-S and gives corresponding optimal feedback control law at a series of time points. Besides, the consistency of estimated MF-S is given.

\subsection{Assumptions and Scenario Descriptions}
$\mathcal{A}_i$ has time-equidistant observations $(l\delta t, l=0,...,N_t)$ on its private trajectory, $T=N_t\delta t$. Agents can change their strategies at time points $0=k_0\delta t<...<k_d\delta t$. At time $k_j\delta t$, $\mathcal{A}_i$ estimates MF-S as $z_i^{k_j}(k_j\delta t)$, with private estimation error $E_i^{k_j}$. 

\emph{A5}: 
At time $k_j\delta t$, $\mathcal{A}_i$ takes $z_i^{k_j}(k_j\delta t)$ as the current actual MF-S, and assumes all agents have the same correct estimation to give its modified feedback control law $u_i(t)=\phi_i^{k_j}(x_i(t),t),k_j\delta t\leq t\leq T$.
\subsection{$\mathcal{A}_i$'s Strategy }
At time $k_j\delta t$, According to \emph{A5} and (\ref{deqn_ex11}), $\mathcal{A}_i$ predicts MF according to 
\begin{equation}
\label{deqn_ex32}
\begin{split}
dz_i^{k_j}(t)=&[(A+C-(B+F)R^{-1}B^TP_0(t))z_i^{k_j}(t)-(B\\
&+F)R^{-1}B^T\mathcal{G}(t)]dt\\
z_i^{k_j}(k_j\delta t)=&z_i^{k_j}(k_j\delta t),\\
\bar{u}_i^{k_j}(t)=&-R^{-1}B^T(P_0(t)z_i^{k_j}(t)+\mathcal{G}(t)).
\end{split}
\end{equation} 
Then its modified feedback control law is 
\begin{equation}
\label{deqn_ex33}
\begin{split}
&u_i(t)=\phi_i^{k_j}(x_i(t),t), k_j\delta t\leq t\leq T\\
&\phi_i^{k_j}(x_i(t),t)=-R^{-1}B^T(P_1(t)x_i(t)+g_i^{k_j}(t)).\\
\end{split}
\end{equation} 
where
\begin{equation}
\label{deqn_ex34}
\begin{split}
dg_i^{k_j}(t)=&-[(A^T-P_1(t)BR^{-1}B^T)g_i^{k_j}(t)+(P_1(t)C\\
&-Q\Gamma)z_i^{k_j}(t)+P_1(t)F\bar{u}_i^{k_j}(t)-Q_Is-Q\eta]dt,\\
g_i^{k_j}(T)=&-\bar{Q}_I\bar{s}-\bar{Q}(\bar{\Gamma}z_i^{k_j}(T)+\bar{\eta}).
\end{split}
\end{equation}
\subsection{Estimation Error Affections}
At time $k_j\delta t,j=0,..,d,k_{d+1}=N_t$, the private error of $\mathcal{A}_i$ is $E_i^{k_j}$, the mean error is $\bar{E}_{k_j}$, define $E_{k_j}^i:=(\bar{E}_{k_j}^T,(E_i^{k_j})^T)^T$. Define $z_{k_j},\bar{u}_{k_j},g_{k_j}$ as $z_i^{k_j},\bar{u}_i^{k_j},g_i^{k_j}$ under complete information. 

Define $\Delta z_i^{k_j}(t):=z_i^{k_j}(t)-z_{k_j}(t)$, $\Delta g_i^{k_j}(t):=g_i^{k_j}(t)-g_{k_j}(t)$, $\Delta z_A^{k_j}(t):=z_A(t)-z_{k_j}(t)$. Then according to (\ref{deqn_ex13}), (\ref{deqn_ex16}) and (\ref{deqn_ex18}), we have
\begin{equation}
\label{deqn_ex35}
\begin{split}
\Delta z_i^{k_j}(t)=&\Phi_1(t)\Phi_1^{-1}(k_j\delta t)E_i^{k_j},\\
\Delta g_i^{k_j}(t)=&\mathcal{M}_g^{k_j}(t)E_i^{k_j},\\
\Delta z_A^{k_j}(t)=&\mathcal{M}_z^{k_j}(t)\bar{E}_{k_j},
\end{split}
\end{equation}
where 
\begin{equation}
\label{deqn_ex36}
\begin{split}
\mathcal{M}_g^{k_j}(t)=&\mathcal{M}_g(t)\Phi_1(0)\Phi_1^{-1}(k_j\delta t),\\
\mathcal{M}_z^{k_j}(t)=&-\Phi_z(t)\int_{k_j\delta t}^t\Phi_z^{-1}(s)(B+F)R^{-1}B^T\mathcal{M}_g^{k_j}(s)ds.\\
\end{split}
\end{equation}
The private trajectory during $t\in [k_j\delta t,k_{j+1}\delta t]$ satisfies
\begin{equation}
\label{deqn_ex37}
\begin{split}
dx_i(t)=&[(A-BR^{-1}B^TP_1(t))x_i(t)+\mathcal{B}_i^{k_j}(t)]dt\\
&+\mathcal{K}^{k_j}(t)E_{k_j}dt+DdW_i(t),\\
\end{split}
\end{equation}
where
\begin{equation}
\label{deqn_ex38}
\begin{split}
\mathcal{K}_i^{k_j}(t)=&\mathcal{K}_i(t)\Phi_1(0)\Phi_1^{-1}(k_j\delta t),\\
\bar{\mathcal{K}}_{k_j}(t)=&(C-FR^{-1}B^TP_1(t))\mathcal{M}_z^{k_j}(t)-FR^{-1}B^T\mathcal{M}_g^{k_j}(t),\\
\mathcal{B}_i^{k_j}(t)=&(C-FR^{-1}B^TP_1)z_i^{k_j}(t)-(B+F)R^{-1}B^Tg_i^{k_j}(t).\\
\end{split}
\end{equation} 
\subsection{Error-Based State Estimation Method}
At time $k_{j+1}\delta t$, $j=0,...,d-1$, $\mathcal{A}_i$ gets discrete observations $x_i((k_j+l)\delta t),l=0,...,k_{j+1}-k_j$. First, $\mathcal{A}_i$ estimates $E_{k_j}^i$ according to the error estimation method proposed in the last section. The transition probability density of $x_i$ during $t\in [k_{j}\delta t,k_{j+1}\delta t]$ is
\begin{equation}
\label{deqn_ex39}
\begin{split}
\mathcal{T}_{k_j}(t,y|s,x;E_{k_j}^i)=\frac{1}{\sqrt{(2\pi)^d|\Sigma|}}\exp^{-\frac{1}{2}(y-\mu_{k_j})^T\Sigma^{-1}(y-\mu_{k_j})},
\end{split}
\end{equation}
where
\begin{equation}
\label{deqn_ex40}
\begin{split}
&\mu_{k_j}(t|s,x;E_{k_j}^i)=F_{k_j}(t|s,x)+\Phi_x(t)\int_s^t\Phi_x^{-1}\mathcal{K}_{k_j}d\omega E_{k_j}^i,\\
&F_{k_j}(t|s,x)=\Phi_x(t)\Phi_x^{-1}(s)x+\Phi_x(t)\int_s^t\Phi_x^{-1}(\omega)\mathcal{B}_i^{k_j}(\omega)d\omega.
\end{split}
\end{equation}
Represent $x_i(s\delta t), s\in\{0,1,...,N_t\}$ as $x_i^s$, the likelihood function is
\begin{equation}
\label{deqn_ex41}
\begin{split}
\mathcal{L}_{k_j}(E^i_{k_j}|x_i)=\prod_{l=k_j+1}^{k_{j+1}}\mathcal{T}_{k_j}(l\delta t,x_i^l|(l-1)\delta t,x_i^{l-1};E_{k_j}^i).
\end{split}
\end{equation}
and the log-likelihood function
\begin{equation}
\label{deqn_ex42}
\begin{split}
l_{k_j}(E_{k_j}^i|x_i)=\sum_{l=k_j+1}^{k_{j+1}}log(\mathcal{T}_{k_j}(l\delta t,x_i^l|(l-1)\delta t,x_i^{l-1};E_{k_j}^i)),
\end{split}
\end{equation}
then $\mathcal{A}_i$'s maximum likelihood estimation of $E_{k_j}^i$ is
\begin{equation}
\label{deqn_ex43}
\begin{split}
&\hat{E}_{k_j}^i=\arg\max_{E^i_{k_j}\in\Lambda}l_{k_j}(E_{k_j}^i|x_i).\\
\end{split}
\end{equation}

Then, for estimated $\hat{E}_i^{k_j},\hat{\bar{E}}^i_{k_j}$, $\mathcal{A}_i$ can estimate $z_A(t),k_j\delta t\leq k_{j+1}\delta t$ by
\begin{equation}
\label{deqn_ex44}
\begin{split}
z_i^{k_{j+1}}(t)=&z_i^{k_j}(t)+\mathcal{M}_z^{k_j}(t)\hat{\bar{E}}^i_{k_j}-\Phi_1(t)\Phi_1(k_j\delta t)^{-1}\hat{E}_i^{k_j}.
\end{split}
\end{equation}

\subsection{Initial Error Tolerant Distributed Mean Field Control}
Combining $\mathcal{A}_i$'s state estimation method and feedback control law, the initial error tolerant distributed mean field control (IET-DMFC) algorithm can be established.
 \begin{algorithm}[!ht]
    \renewcommand{\algorithmicrequire}{\textbf{Input: }}
	\renewcommand{\algorithmicensure}{\textbf{Output:}}
	\caption{$\mathcal{A}_i$'s Strategy in IET-DMFC}
    \label{power}
    \begin{algorithmic}[1] 
        \REQUIRE  $z_i(0), x_i(0),\Theta, k_j,j=1,..,d$; 
        \STATE Solve (\ref{deqn_ex7}), (\ref{deqn_ex8}), (\ref{deqn_ex4}), (\ref{deqn_ex14}), (\ref{deqn_ex17}), (\ref{deqn_ex19}), (\ref{deqn_ex27});
        \STATE Predict MF-S $(z_i^{0}(t))_{0\leq t\leq k_{1}\delta t}$ by (\ref{deqn_ex11});
        \STATE Get $(g_i^{0}(t))_{0\leq t\leq k_{1}\delta t}$ by (\ref{deqn_ex5});
        \STATE $u_i(t)\leftarrow\phi_i^{0}(x_i(t),t), 0\leq t\leq k_{1}\delta t$;
        \STATE Evolve according to (\ref{deqn_ex1}), $0\leq t\leq k_{1}\delta t$;
        \FOR{$j=0$ to $d-1$}
            \STATE Get observations $x_i^l,l=k_j,...,k_{j+1}$; 
            \STATE Estimate error $\hat{E}^i_{k_j}$ by (\ref{deqn_ex43});
            \STATE Estimate current MF-S $z_i^{k+1}(k_{j+1}\delta t)$ by (\ref{deqn_ex44});
            \STATE Predict MF-S $(z_i^{k+1}(t))_{k_{j+1}\delta t\leq t\leq T}$ by (\ref{deqn_ex32});
            \STATE Get $(g_i^{k_j}(t))_{k_{j+1}\delta t\leq t\leq k_{j+2}\delta t}$ by (\ref{deqn_ex34});
            \STATE $u_i(t)\leftarrow\phi_i^{k_{j+1}}(x_i(t),t), k_{j+1}\delta t\leq t\leq k_{j+2}\delta t$;
            \STATE Evolve according to (\ref{deqn_ex1}), $k_{j+1}\delta t\leq t\leq k_{j+2}\delta t$;
        \ENDFOR
    \end{algorithmic}
\end{algorithm}  
\subsection{Consistency of the Estimation for MF-S}
The following theorem gives the consistent property of the state estimation method.

{\bf{Theorem 4.1}} When $A1-A5$ holds, define $\kappa_j:=k_j-k_{j-1}$, estimated MF-S $z_i^{k_{j}}(k_{j}\delta t), j=1,...,d$ satisfies
\begin{equation}
\label{deqn_ex45}
\begin{split}
z_i^{k_{j}}(k_{j}\delta t)\stackrel{a.s.}{\longrightarrow} z_A(k_{j}\delta t) \ as\ \kappa_j\rightarrow \infty
\end{split}
\end{equation}
Proof: according to Lemma 3.1, $\exists\Omega_0\in \mathcal{F}, P(\Omega_0)=1, \forall \omega \in \Omega_0$, 
\begin{equation*}
\lim_{\kappa_j\rightarrow \infty}|\hat{E}^i_{k_{j-1}}(\omega)-E^i_{k_{j-1}}|=0.
\end{equation*}
Define $\delta z_i^j:=z_i^{k_{j}}(k_{j}\delta t)-z_A(k_{j}\delta t)$, $\delta E_{k_j}^i:=\hat{E}^i_{k_{j-1}}-E^i_{k_{j-1}}$, according to (\ref{deqn_ex44}), 
\begin{equation*}
\delta z_i^j=[\mathcal{M}_z^{k_{j-1}}(k_{j}\delta t),-\Phi_1(k_{j}\delta t)\Phi_1(k_{j-1}\delta t)^{-1}]\delta E_{k_j}^i.
\end{equation*}
Applying the continuity of $\mathcal{M}_g(t)$,$\Phi_1(t),\Phi_z(t)$ on $t\in[0,T]$, for fixed $T$, the continuity and boundedness of $\mathcal{M}_g^{k_{j-1}}$ and $\mathcal{M}_z^{k_{j-1}}$ can be derived. Therefore $\exists C>0$,
\begin{equation*}
|\mathcal{M}_z^{k_{j-1}}(k_{j}\delta t),-\Phi_1(k_{j}\delta t)\Phi_1(k_{j-1}\delta t)^{-1}|<C, \forall j.
\end{equation*}
Then we have $|\delta z_i^j(\omega)|<C|\delta E_{k_j}^i(\omega)|$, and
\begin{equation*}
\lim_{\kappa_j\rightarrow \infty}|\delta z_i^j(\omega)|=0.
\end{equation*}
Which is,
\begin{equation*}
z_i^{k_{j}}(k_{j}\delta t)\stackrel{a.s.}{\longrightarrow} z_A(k_{j}\delta t) \ as\ \kappa_j\rightarrow \infty
\end{equation*}

$\Box$

Above theorem shows that, when observation frequency tends to infinity, the estimation error $E^i_{k_j}$ at time $k_j\delta t$ tends to $0$. So \emph{A5} is a reasonable and self-consistent assumption in this sense.
\section{Simulations}
We select $100$ agents from the population. Set $T=2,\delta t=0.02, A=R=s=\bar{s}=1,C=-1,B=F=Q_I=\bar{Q}_I=0.5,Q=\bar{Q}=\eta=\bar{\eta}=0.1,\Gamma=\bar{\Gamma}=1$. The initial distribution is a normal distribution with $z^0=0$ as the expectation and $0.1$ as the variance. Initial private errors $\{E_i,1\leq i\leq 100\}$ conform to a normal distribution with $\bar{E}=10$ as the expectation and $2$ as the variance.
\begin{figure}[!t]
\centering
\subfigure[]{\label{trajc}\includegraphics[width=1.65in]{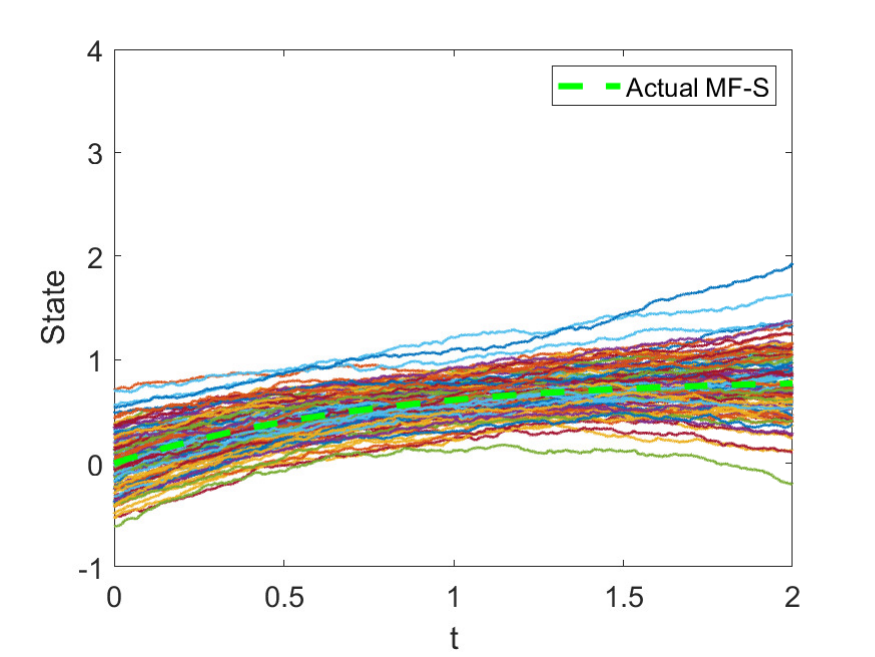}}
\hfill 
\subfigure[]{\label{trajA}\includegraphics[width=1.65in]{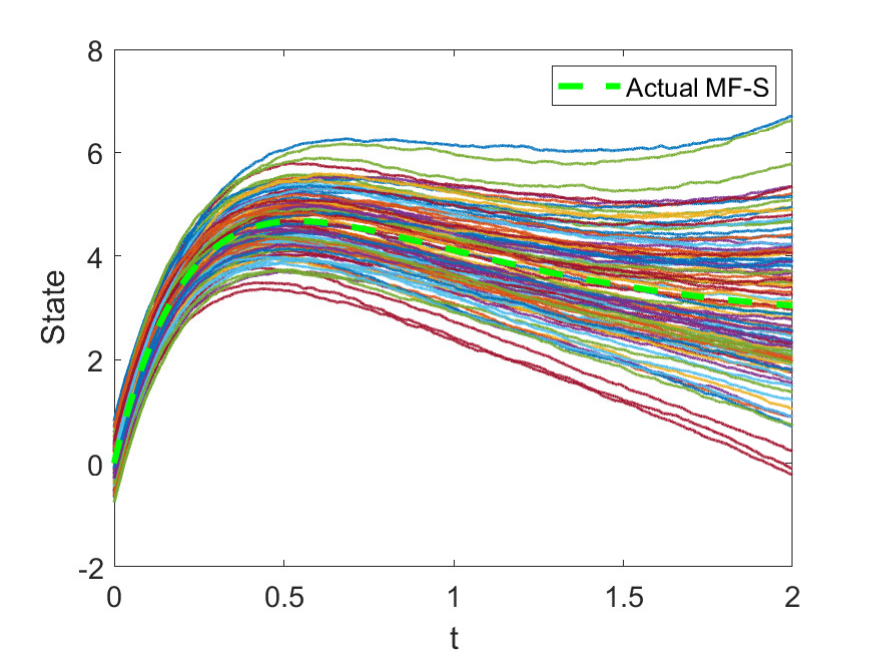}}

\subfigure[]{\label{trajm1}\includegraphics[width=1.65in]{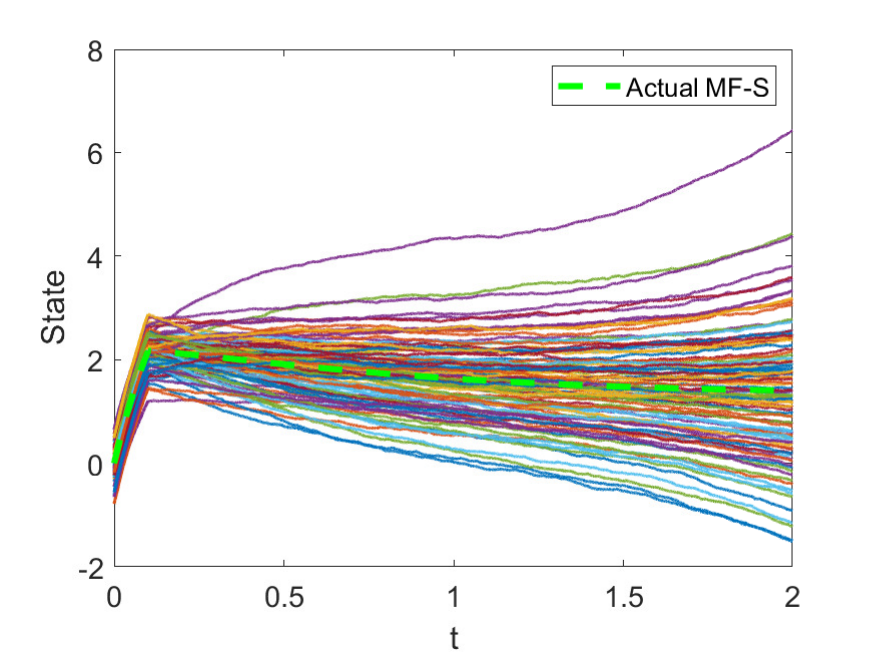}}
\hfill 
\subfigure[]{\label{MFScompare}\includegraphics[width=1.65in]{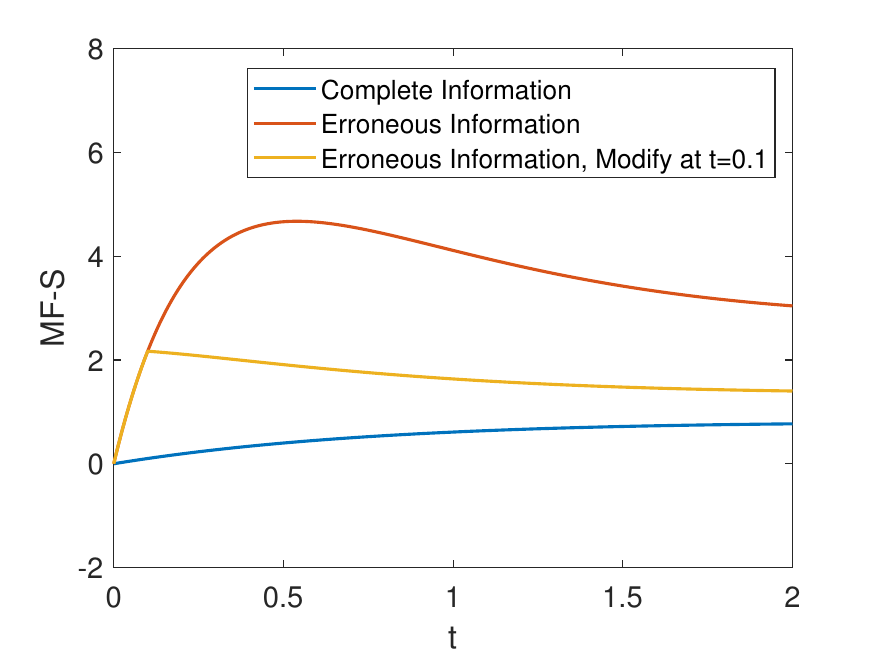}}
\caption{Agents' trajectories under different settings. (a) Complete information. (b) Erroneous information. (c) Erroneous information, agents take strategy modifications at $t=0.1$. (d) Mean field states under different settings.}
\label{fig_1}
\end{figure}

Fig.\ref{fig_1} illustrates agents' trajectories under different settings. It can be seen that the initial information errors evert a profound impact on population dynamics. In Fig.\ref{trajm1}, agents estimate current MF-S and modify their strategies at $t_1=0.1$. The modified trajectories are more similar to the trajectories under complete information than the trajectories under erroneous information without strategy modification, which is also shown in Fig.\ref{MFScompare}, where the MF-S under different situations are compared. 

\begin{figure}[!t]
\centering
\subfigure[]{\label{meanerror}\includegraphics[width=1.65in]{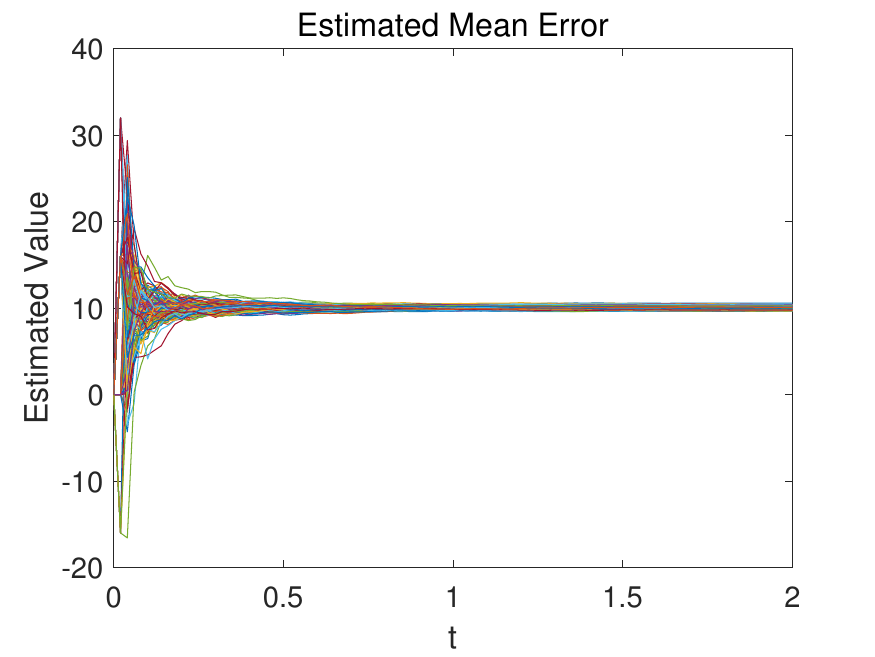}}
\hfill 
\subfigure[]{\label{privateerror}\includegraphics[width=1.65in]{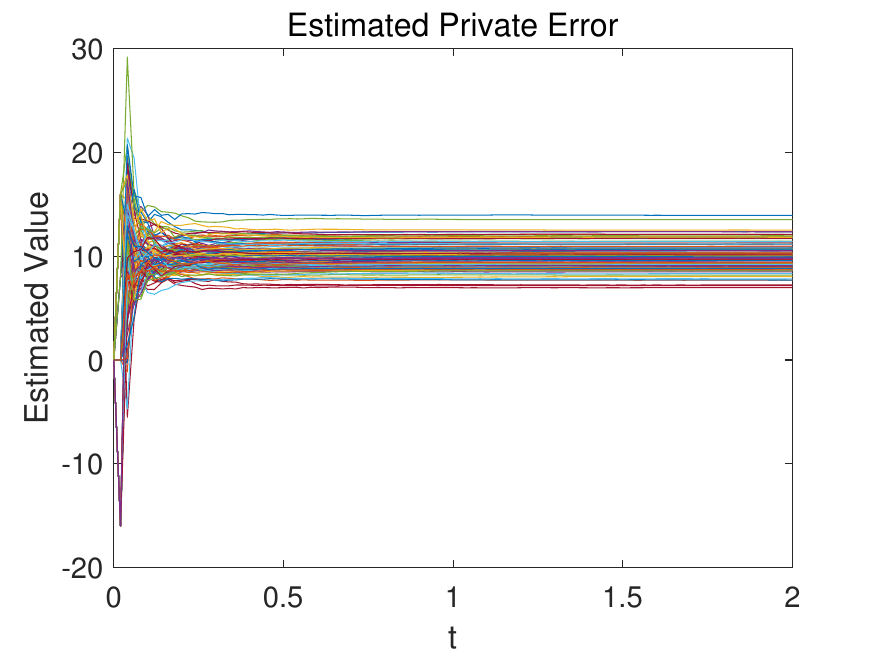}}

\caption{Estimations for initial information errors. (a) Estimations for mean error. (b) Estimations for private errors.}
\label{fig_2}
\end{figure}

Fig.\ref{fig_2} shows the consistent property of our error estimation method. At $t_j=0.02j, j=1,...,100$, $\mathcal{A}_i$ computes the estimation for initial information errors using current discrete observations $x_i(0.02l),l=1,...,j$. In Fig.\ref{meanerror}, as the number of observations increases, estimated values for mean error and private errors intend to be stable and consistent with the true value.

\begin{figure}[!t]
\centering
\subfigure[]{\label{MFP0}\includegraphics[width=1.65in]{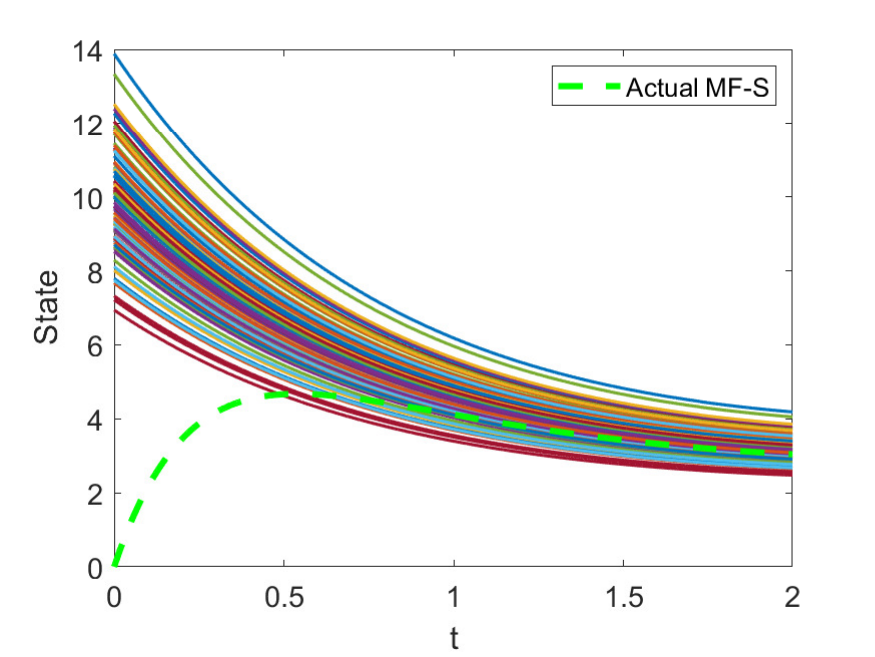}}
\hfill 
\subfigure[]{\label{MFP1}\includegraphics[width=1.65in]{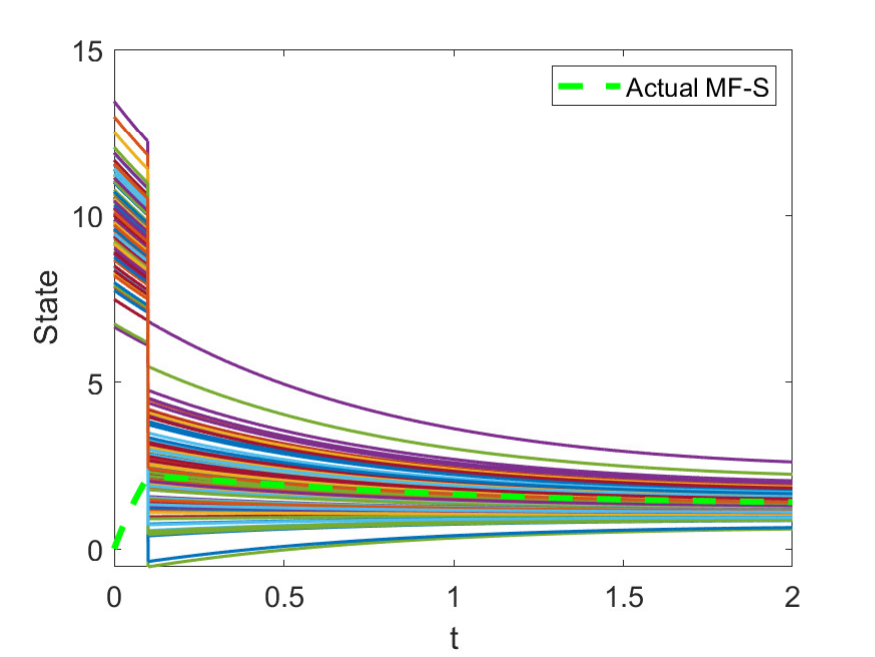}}

\subfigure[]{\label{MFP2}\includegraphics[width=1.65in]{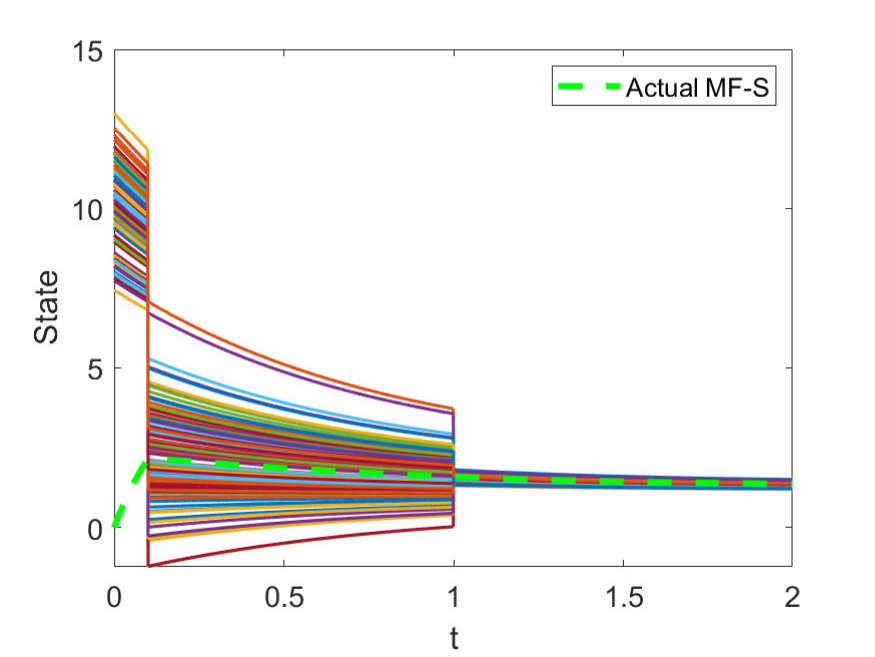}}
\hfill 
\subfigure[]{\label{Stateerror}\includegraphics[width=1.65in]{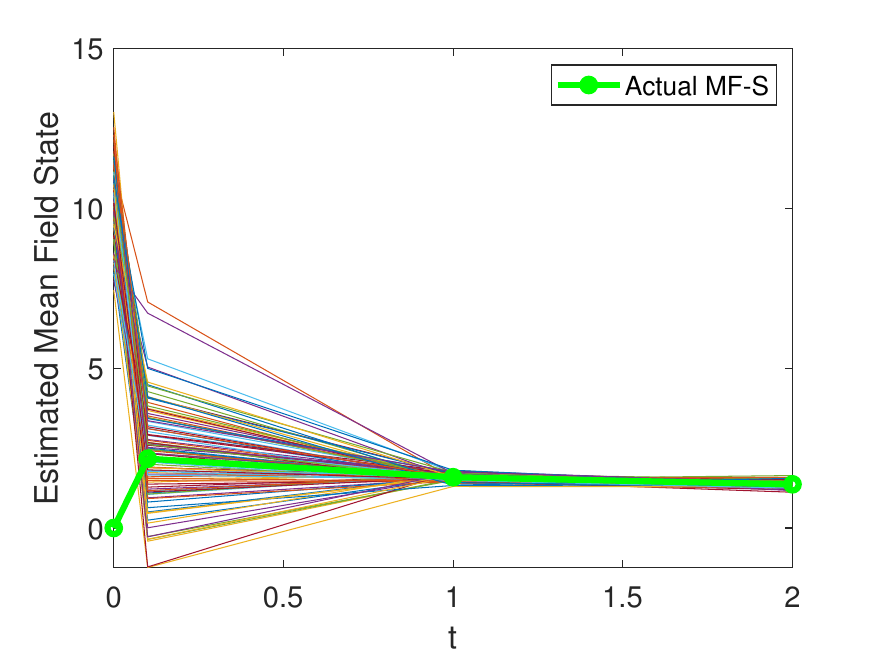}}
\caption{Predicted mean field states under different modification time points settings. (a) Predicted MF-S under erroneous information without strategy modification. (b) Predicted MF-S under erroneous information with strategy modification on $t_1=0.1$. (c) Predicted MF-S under erroneous information with strategy modifications on $t_1=0.1, t_2=1$. (d) Estimations for MF-S in Fig.\ref{MFP2}.}
\label{fig_3}
\end{figure}

We conduct simulations under different modification time points settings. Fig.\ref{fig_3} compare predicted MF-S under different settings. In Fig.\ref{MFP0}, agents predict MF-S according to the initial information, and heterogeneous initial errors cause heterogeneous predictions for MF. In Fig.\ref{MFP1}, agents estimate current MF-S and predict future MF-S at time $t_1=0.1$. Estimations are closer to the actual MF-S, but due to the small amount of observations, estimations still differ significantly among individuals. In Fig.\ref{MFP2} and Fig.\ref{Stateerror}, agents modify their strategies at $t_1=0.1,t_2=1$. Compared with the estimations at $t_1$, the estimations and predictions converges to the actual value significantly at $t_2$. 

Above simulations demonstrate the effectiveness of IET-DMFG in information error tolerance and strategy modifications, and verify the consistency conclusion we mentioned in Lemma 3.1 and Theorem 4.1.

\section{Conclusion}
This paper proposes an initial error tolerant distributed mean field control (IET-DMFC) framework for systems with partial and discrete observations. By analyzing error propagation in LQMFGs, we develop a distributed MLE-based error estimation method and a segmented state estimation method, enabling each agent to infer and correct initial errors using only discrete information of its private trajectory. For future work, more properties of the estimations such as asymptotic normality can be discussed, and the case of a limited number of agents can be considered. 

\end{document}